\documentclass[11pt]{article}  

\usepackage{amsfonts}

\begin{document}

\title{Bivariate Daubechies Scaling Functions\footnote{This work was partially supported by National Science Foundation grant DMS-9820221.}}        
\author{
Edward Aboufadel\footnote{Grand Valley State University, e-mail: aboufade@gvsu.edu},
Amanda Cox\footnote{St. Olaf College, e-mail: coxa@stolaf.edu},
and Amy Vander Zee\footnote{Grand Valley State University, e-mail: vanderza@river.it.gvsu.edu }}
\date{Technical Report:  \today \quad draft}            
\maketitle

\section{Introduction}       

\subsection{Background about Wavelets}

Wavelets are functions used to approximate other functions or data. They are particulary effective in approximating functions with discontinuites or sharp changes. Wavelets further preserve both the global picture and local detail of a function. Applications of wavelets include compressing images, such as X-rays and the FBI's fingerprint collection, and denoising data, found in many areas, including geology, meteorology, astronomy, acoustics, and econometrics. 

\subsection{Discussion of $D_4$}

$D_4$ is the univariate Daubechies scaling function with four refinement coeffiecients. $D_4$ may be used to exactly reproduce univariate functions which are constant or linear.

\section{Overview of Bivariate Daubechies Scaling Function}

\subsection{Conditions on $\phi(x,y)$}

The idea of this project is to create scaling functions $\phi(x,y)$ of two variables which satisfy properties of the type satisfied by $D_4$.  The \emph{dilation equation} for a bivariate Daubechies scaling function {\bf [D]} will be:
\begin{equation} \label{E:Dilation}
\phi(x,y) = \sum_i\sum_j c_{i,j}\phi(2x-i,2y-j).
\end{equation}
\noindent  where the sums are over all integers.  This scaling function will satisfy several conditions.  One condition to satisfy is \emph {compact support}.  The support of the bivariate scaling function is from $0$ to $3$ in the $x$ direction and from $0$ to $3$ in the $y$ direction.  The function $\phi$ is zero outside of this support, and is also zero on the boundary.  There is also an \emph{averaging condition} that comes from the idea that the average value of the scaling function is not 0.  Related to this condition is the assumption that the inner product
\begin{equation} \label{E:ScalNorm}
<\phi(x,y),\phi(x,y)> = \int_{-\infty}^\infty\int_{-\infty}^\infty\phi(x,y)^2 dx dy = 1
\end{equation}

A bivariate scaling function must be orthogonal to its three types of integer translates: those in the $x$ direction only, those in the $y$ direction only, and those in both the $x$ and $y$ directions.  This is called the \emph{orthogonality condition}.  (The orthogonality is important in order that the function and its translates form a basis.)  Also, the scaling function must be designed to exactly reproduce linear functions of the form $f(x,y) = kx + ly + m$, where $k$, $l$, and $m$ are constants.  This last condition is called the \emph{regularity condition}, and it also guarantees that $\phi$ is continuous.  (However, Han points out that any $\phi$ we create will not be differentiable everywhere, a property shared by $D_4$ {\bf [H]}.) The consequences of all four conditions will be explored in the next several sections.

\subsection{A Formula for the Refinement Coefficients}

Due to the orthogonality condition, the set of functions
\begin{equation}
\{ \phi(2x-i, 2y-j) | i, j \in {\mathbb Z} \}
\end{equation}
\noindent is a basis of a vector space of which $\phi(x,y)$ is a member.  This observations allows us to prove this Lemma:

{\bf Lemma 1:} For all integers $i$ and $j$, $c_{i,j}=4<\phi(x,y),\phi(2x-i,2y-j)>$.

{\bf Proof of Lemma 1:}  According to the Orthogonal Decomposition Theorem,  
\begin{equation} \label{E:Decompose}
\phi(x,y)=\sum_i\sum_j \frac{<\phi(x,y),\phi(2x-i,2y-j)>}{<\phi(2x-i,2y-j),\phi(2x-i,2y-j)>}\phi(2x-i,2y-j)
\end{equation}
\noindent and that
\begin{equation} \label{E:Refine}
c_{i,j} =\frac{<\phi(x,y),\phi(2x-i,2y-j)>}{<\phi(2x-i,2y-j)\phi(2x-i,2y-j)>}.
\end{equation}

The denominator in this formula turns out to equal $\frac{1}{4}$.  To see this, observe that the definition of inner product, along with (\ref{E:ScalNorm}) yields
\begin{eqnarray*}
<\phi(2x-i,2y-j)\phi(2x-i,2y-j)> 
&=& \int_{-\infty}^\infty \int_{-\infty}^\infty \phi(2x-i,2y-j)\phi(2x-i,2y-j) dx dy \\
&=& \frac{1}{4}\int_{-\infty}^\infty \int_{-\infty}^\infty \phi(u-i,v-j)\phi(u-i,v-j) du dv \\
&=& \frac{1}{4}<\phi(x,y),\phi(x,y)>.
\end{eqnarray*}
\noindent Therefore, 
\begin{equation} \label{E:Because}
<\phi(2x-i,2y-j)\phi(2x-i,2y-j)> = \frac{1}{4},
\end{equation}
\noindent and we have proven that $c_{i,j}=4<\phi(x,y),\phi(2x-i,2y-j)>$.

\subsection{Consequence of Compact Support:  Only 16 Non-Zero Refinement Coefficients}

Next, assuming that the compact support for a bivariate Daubechies scaling function is $0\leq x\leq3$ and $0\leq y\leq 3$, we prove this lemma:

{\bf Lemma 2:} For $i \not= 0,1,2,3$ or $j \not= 0,1,2,3$, $c_{i,j}=0$.

{\bf Proof of Lemma 2:}  From the definition of inner product, we know that
$$<\phi(x,y),\phi(2x-i,2y-j)>=\int_{-\infty}^\infty \int_{-\infty}^\infty \phi(x,y)\phi(2x-i,2y-j) dx dy.$$
\noindent When $i\leq -2$ or $i\geq 6$, and $j\leq -2$ or $j\geq 6$, then $\phi(2x-i,2y-j)$ will be non-zero where $\phi(x,y)=0$, and vice versa.  Therefore, we can definitely conclude that $c_{i,j} = 0$ unless $-2\leq i \leq5$ and $-2\leq j\leq5$.

It remains to show that $c_{i,j} = 0$ for $i = -2, -1, 4,$ or 5, and $j = -2, -1, 4,$ or 5.

We can use our knowledge that some refinement coefficients are zero in order to prove that other ones are zero.  The key to showing that the other coefficients are zero is the following observation: for any integers $m$ and $n$,
\begin{equation} \label{E:RefSquare}
\phi(\frac{m}{2}, \frac{n}{2}) = c_{m-1, n-1}\phi(1,1) +
c_{m-1, n-2}\phi(1,2) + c_{m-2, n-1}\phi(2,1) + c_{m-2, n-2}\phi(2,2)
\end{equation}
\noindent due to the compact support requirement.  For instance,
\begin{equation} \label{E:Elim22}
\phi(-\frac{1}{2}, -\frac{1}{2}) = c_{-2, -2}\phi(1,1) +
c_{-2, -3}\phi(1,2) + c_{-2, -3}\phi(2,1) + c_{-3, -3}\phi(2,2).
\end{equation}
\noindent  Since any refinement coefficient which has -3 as one of its indices is already known to be zero, and due to the compact support condition, (\ref{E:Elim22}) becomes
\begin{equation}
0 = c_{-2, -2}\phi(1,1).
\end{equation}
Already we can conclude that $c_{-2, -2} = 0$ or $\phi(1,1) = 0$.

From an assumption about one of $\phi(1,1)$, $\phi(1,2)$, $\phi(2,1)$ or $\phi(2,2)$ being non-zero, we can prove that more than half of the remaining refinement coefficients that we are trying to eliminate are, in fact, zero.  We will demonstate how this works if we assume $\phi(1,1) \not= 0$.

If $\phi(1,1) \not= 0$, then $c_{-2, -2} = 0$.  Then using $m = -1$ and $n=0$, (\ref{E:RefSquare}) becomes
\begin{equation}
\phi(-\frac{1}{2}, 0) = c_{-2, -1}\phi(1,1) +
c_{-2, -2}\phi(1,2) + c_{-3, -1}\phi(2,1) + c_{-3, -2}\phi(2,2). 
\end{equation}
\noindent  So, by previous results, $0 = c_{-2, -1}\phi(1,1)$, and therefore $c_{-2,-1} = 0.$
In a similar manner, it can be shown that $c_{-2, j} = 0$ for all remaining $j$, and nearly identical arguments can be made to show that $c_{i, -2} = 0$ for all remaining $i$.

Refinement coefficients of the form $c_{i, -1}$ and $c_{-1, j}$ can be disposed of in the same way.  For example, using $m = 0$ and $n = -1$, (\ref{E:RefSquare}) becomes
\begin{equation} 
\phi(0, -\frac{n}{2}) = c_{-1, -2}\phi(1,1) +
c_{-1, -3}\phi(1,2) + c_{-2, -2}\phi(2,1) + c_{-2, -3}\phi(2,2)
\end{equation}  Again, by previous results, this equation reduces to $0 = c_{-1, -2}\phi(1,1)$, so $c_{-1, -2} = 0$.  Other coefficients are then dealt with in the same way.  Consequently, we can eliminate more than half of the refinement coefficients that we are trying to show are zero by an assumption that $\phi(1,1) \not= 0$.  The same can be said of the other three values of $\phi$.

In the case where both $\phi(1,1)$ and $\phi(2,2)$ are non-zero, then all of the unwanted coefficients are clearly eliminated.  The same is true if both $\phi(1,2)$ and $\phi(2,1)$ are non-zero.

Without loss of generality, if $\phi(1,1)$ and $\phi(2,1)$ are both zero, and the other two are non-zero, then a few more steps are required to reach our result.  In this case, there are eight refinement coefficients that are not determined to be zero by our previous arguments:  $c_{i, j}$ where $i$ can be 0, 1, 2, or 3, and $j$ is $-2$ or $-1$.  However, the same type of argument can be used to eliminate these one at a time.  For instance, using $m = 1$ and $n = 0$, (\ref{E:RefSquare}) becomes
\begin{equation} 
\phi(\frac{1}{2}, 0) = c_{0, -1}\phi(1,1) +
c_{0, -2}\phi(1,2) + c_{-1, -1}\phi(2,1) + c_{-1, -2}\phi(2,2).
\end{equation}  The first and third terms are zero by our assumption in this case.  The coefficient $c_{-1, -2}$ is zero follows from the assumption that $\phi(1,2)$ is non-zero, using our arguments above.  Therefore, we have that $0 = c_{0, -2}\phi(1,2)$, and we can conclude that $c_{0, -2} = 0$.  All eight can be eliminated in this way.

Finally, if three of the four function values are zero, we can use the fourth value to eliminate more than half of the coefficients and then arguments such as in the previous paragraph to show that the rest are zero.  We omit the details.

If all four of the function values are zero, then $\phi$ is the zero function.

In summary, for all cases, all of the refinement coefficients are zero unless $0 \leq i \leq 3$ and $0 \leq j \leq 3$.

\section{Deriving the 14 equations}

\subsection{The Averaging Equation}

In this section we prove that
\begin{equation} \label{E:Averaging}
\sum_i\sum_j c_{i,j}=4.
\end{equation}

Here we will use the averaging assumption, which states that
\begin{equation} \label{E:Averaging2}
\int_{-\infty}^\infty \int_{-\infty}^\infty \phi(x,y) dx dy =1.
\end{equation}
\noindent Combining this with (\ref{E:Dilation}) and (\ref{E:Because}), we can calculate the following:

\begin{eqnarray*}
\int_{-\infty}^\infty \int_{-\infty}^\infty \phi(x,y)dx dy &=&\int_{-\infty}^\infty \int_{-\infty}^\infty(\sum_i\sum_j c_{i,j}\phi(2x-i,2y-j)) dx dy\\
&=&\frac{1}{4}\int_{-\infty}^\infty \int_{-\infty}^\infty(\sum_i\sum_j c_{i,j} \phi(x,y))dx dy \\
&=&\frac{1}{4}\sum_i\sum_j c_{i,j}\int_{-\infty}^\infty \int_{-\infty}^\infty \phi(x,y)dx dy
\end{eqnarray*}
\noindent  Note that the last step is valid because the sums are finite.  As a result we have the following equation:
$$\int_{-\infty}^\infty \int_{-\infty}^\infty \phi(x,y)dx dy = \frac{1}{4} \sum_i\sum_j c_{i,j}\int_{-\infty}^\infty \int_{-\infty}^\infty \phi(x,y)dx dy$$
\noindent  Due to our averaging condition, we can divide both sides by the double integral, leaving us with $1 =\frac{1}{4}\sum_i\sum_j c_{i,j}$, or $4 =\sum_i\sum_j c_{i,j}$.

\subsection{The Orthogonality Equations I}

By the orthogonality condition, we assume that translates in the $x$ or $y$ direction (or both) are orthogonal to the scaling function.  Using $b$ and $d$ to translate the variables, the dilation equation (\ref{E:Dilation}) becomes:
$$\phi(x-b,y-d) = \sum_i\sum_jc_{i,j}\phi(2(x-b)-i, 2(y-d)-j).$$
\noindent Letting $m=2b+i$ and $n=2d+j$ leads to
$$\phi(x-b,y-d) = \sum_{m}\sum_{n}c_{m-2b,n-2d}\phi(2x-m, 2y-n).$$

We are interested in the inner product of $\phi(x,y)$ and $\phi(x-b,y-d)$. That is,
\begin{eqnarray*}
&<& \phi(x,y),\phi(x-b,y-d)> \\
&=&<\sum_i\sum_j c_{i,j}\phi(2x-i,2y-j),\sum_m\sum_nc_{m-2b,n-2d}\phi(2x-m,2y-n)> \\
&=&\sum_i \sum_j \sum_m \sum_n c_{i,j}c_{m-2b,n-2d}<\phi(2x-i,2y-j),\phi(2x-m,2y-n)> \\
&=& \frac{1}{4} \sum_i \sum_j \sum_m \sum_n c_{i,j}c_{m-2b,n-2d}<\phi(x-i,y-j),\phi(x-m,y-n)> \\
\end{eqnarray*}
\noindent with the last line being true due to (\ref{E:Because}).  

When $i\not=m$ or $j\not=n$, $<\phi(x-i,y-j),\phi(x-m,y-n)>=0$, by the orthogonality condition.  So, in the previous paragraph, $\sum_m \sum_n$ only counts when $i=m$ and $j=n$, and it is actually the case that $<\phi(x-i,y-j),\phi(x-m,y-n)>=1$.  So we have that
\begin{equation}
<\phi(x,y),\phi(x-b,y-d)> = \sum_i \sum_j c_{i,j} c_{i-2b,j-2d}.
\end{equation}
\noindent If $b$ and $d$ are not both zero, then
\begin{equation} \label{E:PreOrthog}
0=\sum_i\sum_jc_{i,j}c_{m-2b,n-2d}.
\end{equation}

Specific values of $b$ and $d$ will lead to specific \emph{orthogonality equations}.  These equations take into account Lemma 2, which showed that only 16 refinement coefficients are non-zero.  If $b$ and $d$ are both 1, then (\ref{E:PreOrthog}) becomes
\begin{equation} \label{E:Orthog2} 
c_{3,3}c_{1,1}+c_{3,2}c_{1,0}+c_{2,3}c_{0,1}+c_{2,2}c_{0,0}=0
\end{equation}
\noindent  Similarly, if $b\not=0$ and $d=0,$ or if $b=0$ and $d\not=0$, respectively, (\ref{E:PreOrthog}) reduces to
\begin{equation} \label{E:Orthog3}
c_{3,3}c_{1,3}+c_{3,2}c_{1,2}+c_{2,3}c_{0,3}+c_{2,2}c_{0,2}+
c_{3,1}c_{1,1}+c_{3,0}c_{1,0}+c_{2,1}c_{0,1}+c_{2,0}c_{0,0}= 0
\end{equation}
\noindent and 
\begin{equation} \label{E:Orthog4}
c_{3,3}c_{3,1}+c_{3,2}c_{3,0}+c_{2,3}c_{2,1}+c_{2,2}c_{2,0}+
c_{1,3}c_{1,1}+c_{1,2}c_{1,0}+c_{0,3}c_{0,1}+c_{0,2}c_{0,0}=0
\end{equation}
\noindent  No other orthogonality equations can be derived from (\ref{E:PreOrthog}).

\subsection{The Orthogonality Equations II}

There is one more equation to be derived based on the orthogonality condition.  Applying the dilation equation (\ref{E:Dilation}) to (\ref{E:ScalNorm}) gives us:
\begin{eqnarray*}
<\phi(x,y),\phi(x,y)>
&=&<\sum_i\sum_j c_{i,j}\phi(2x-i, 2y-j),\sum_k\sum_l c_{k,l}\phi(2x-k, 2y-l)> \\
&=&\sum_i\sum_j \sum_k \sum_l c_{i,j} c_{k,l} <\phi(2x-i, 2y-j),\phi(2x-k, 2y-l)> \\
&=&\sum_i\sum_j c_{i,j}^2<\phi(2x-i, 2y-j),\phi(2x-i, 2y-j)>\\
&=&\frac{1}{4}\sum_i\sum_j c_{i,j}^2<\phi(x-i,y-j),\phi(x-i,y-j)> \\ 
\end{eqnarray*}
\noindent with the last two lines being true by the same type of argument that led to (\ref{E:PreOrthog}).  Since 
$$<\phi(x-i,y-j),\phi(x-i,y-j)> = <\phi(x,y),\phi(x,y)>,$$
\noindent we have shown that
\begin{equation} \label{E:AlmostOrthog}
<\phi(x,y),\phi(x,y)>=\frac{1}{4}\sum_i\sum_j c_{i,j}^2<\phi(x,y),\phi(x,y)>
\end{equation}

Due to (\ref{E:ScalNorm}), $<\phi(x,y),\phi(x,y)>$ can be divided out of both sides of (\ref{E:AlmostOrthog}), leaving $\frac{1}{4}\sum_i\sum_j c_{i,j}^2 = 1$, or
\begin{equation} \label{E:Orthog1}
\sum_i\sum_j c_{i,j}^2=4.
\end{equation}
\noindent Due to Lemma 2, the sum in this last equation is made up of 16 terms.

\subsection{The Regularity Equations}

In {\bf [J]} Jia states that, for a positive integer $k$, a sequence $a$ on 
$\mathbb{Z}^s$ satisfies the \emph{sum rules} of order $k$ if $$\sum_{\beta \in \mathbb{Z}^s}\alpha(2\beta + \epsilon)p(2\beta + \epsilon) = \sum_{\beta \in Z^s}\alpha(2\beta)p(2\beta)\ \ \ \ \ \ \forall \epsilon \in \Omega, p \in \Pi_{k-1},$$
where $\Pi_{k-1}$ is the set of polynomials with total degree less than $k$ and $\Omega$ is the set of vertices of the unit cube $[0,1]^s$.  Jia proves that satisfying the sum rules is equivalent to satisfying the regularity condition.  We can apply the result of Jia to the case where $s=2$ (the bivariate case) and $k=2$ (linear functions are reproduced).  Then nine regularity equations follow from the sum rules.

When $p=1$ and $\epsilon=(0,1)$, $\sum c_{2i,2j+1} = \sum c_{2i,2j}$, or
\begin{equation} \label{E:Reg1}
c_{0,1}+c_{0,3}+c_{2,1}+c_{2,3} = c_{0,0}+c_{0,2} + c_{2,0} + c_{2,2}
\end{equation}

When $p=1$ and $\epsilon=(1,1)$, $\sum c_{2i+1,2j+1} = \sum c_{2i,2j}$, or
\begin{equation} \label{E:Reg2}
c_{1,1}+c_{1,3}+c_{3,1}+c_{3,3} = c_{0,0}+c_{0,2} + c_{2,0} + c_{2,2}
\end{equation}

When $p=1$ and $\epsilon=(1,0)$, $\sum c_{2i+1,2j} = \sum c_{2i,2j}$, or
\begin{equation} \label{E:Reg3}
c_{1,0}+c_{1,2}+c_{3,0}+c_{3,2}=c_{0,0}+c_{0,2} + c_{2,0} + c_{2,2}
\end{equation}

When $p=x$ and $\epsilon=(0,1)$, $\sum 2ic_{2i,2j+1} = \sum2ic_{2i,2j}$, or
\begin{equation} \label{E:Reg4}
2c_{2,0}+2c_{2,2} = 2c_{2,1} + 2c_{2,3}
\end{equation}

When $p=x$ and $\epsilon=(1,1)$, $\sum 2i+c_{2i+1,2j+1} = \sum2ic_{2i,2j}$, or
\begin{equation} \label{E:Reg5}
c_{1,1} + c_{1,3} + 3c_{3,1} + 3c_{3,3} = 2c_{2,1} + 2c_{2,3}
\end{equation}

When $p=x$ and $\epsilon=(1,0)$, $\sum2i+c_{2i+1,2j} = \sum2ic_{2i,2j}v$, or
\begin{equation} \label{E:Reg6}
c_{1,0} + c_{1,2} + 3c_{3,0} + 3c_{3,2} = 2c_{2,1} + 2c_{2,3}
\end{equation}

When $p=y$ and $\epsilon=(0,1)$, $\sum2j+c_{2i,2j+1} = \sum2jc_{2i,2j}$, or
\begin{equation} \label{E:Reg7}
c_{0,1}+c_{2,1}+3c_{0,3}+3c_{2,3} = 2c_{0,2}+2c_{2,2}
\end{equation}

When $p=y$ and $\epsilon=(1,1)$, $\sum2j+c_{2i+1,2j+1} = \sum2jc_{2i,2j}$, or
\begin{equation} \label{E:Reg8}
c_{1,1}+c_{3,1}+3c_{1,3}+3c_{3,3}=2c_{0,2}+2c_{2,2}
\end{equation}

When $p=y$ and $\epsilon=(1,0)$, $\sum2jc_{2i+1,2j}=\sum2jc_{2i,2j}$, or
\begin{equation} \label{E:Reg9}
2c_{1,2}+2c_{3,2}=2c_{0,2}+2c_{2,2}
\end{equation}

Note that all of these equations are linear in the refinement coefficients.

\section{Solution of the 14 equations}

\subsection{Method of Solution of the Fourteen Equations}

We now have fourteen equations with sixteen variables.  To review, here are the fourteen equations:

\emph{Averaging equation:}
$$ \sum_i\sum_jc_{i,j}=4 $$

\emph{Orthogonality equations:}
$$ \sum_i\sum_jc_{i,j}^2=4 $$
$$ c_{3,3}c_{1,1}+c_{3,2}c_{1,0}+c_{2,3}c_{0,1}+c_{2,2}c_{0,0}=0 $$
$$ c_{3,3}c_{1,3}+c_{3,2}c_{1,2}+c_{3,1}c_{1,1}+ c_{3,0}c_{1,3}+c_{2,3}c_{0,3}+c_{2,2}c_{0,2} +c_{2,1}c_{0,1}+c_{2,0}c_{0,0}=0 $$  
$$ c_{0,3}c_{0,1}+c_{1,3}c_{1,1}+c_{2,3}c_{2,1}+ c_{3,3}c_{3,1}+c_{0,2}c_{0,0}+c_{1,2}c_{1,0} +c_{2,2}c_{2,0}+c_{3,2}c_{3,0}=0 $$

\emph{Regularity Equations:}
$$c_{0,1}+c_{0,3}+c_{2,1}+c_{2,3}-c_{0,0}-c_{0,2}-c_{2,0}-c_{2,2}=0 $$
$$c_{1,0}+c_{3,0}+c_{1,2}+c_{3,2}-c_{0,0}-c_{0,2}-c_{2,0}-c_{2,2}=0 $$
$$c_{1,1}+c_{1,3}+c_{3,1}+c_{3,3}-c_{0,0}-c_{0,2}-c_{2,0}-c_{2,2}=0 $$
$$2c_{2,0}+2c_{2,2}-2c_{2,1}-2c_{2,3}=0 $$
$$2c_{1,2}+2c_{3,2}-2c_{0,2}-2c_{2,2}=0 $$
$$c_{1,0}+c_{1,2}+3c_{3,0}+3c_{3,2}-2c_{2,0}-2c_{2,2}=0 $$
$$c_{1,1}+c_{1,3}+3c_{3,1}+3c_{3,3}-2c_{2,0}-2c_{2,2}=0 $$
$$c_{0,1}+c_{2,1}+3c_{0,3}+3c_{2,3}-2c_{0,2}-2c_{2,2}=0 $$
$$c_{1,1}+c_{3,1}+3c_{1,3}+3c_{3,3}-2c_{0,2}-2c_{2,2}=0 $$

Due to the complexity of the problem, \emph{Maple} is unable to solve the system of fourteen equations directly.  As a result, we had to take another approach.  We began by 
observing that the regularity equations and the averaging equation form a non-homogeneous system of linear equations.  Solving this system resulted in solutions for the ten following variables $c_{1,2}, c_{2,1}, c_{2,0}, c_{0,2}, c_{0,1}, c_{1,0}, c_{1,1}, c_{3,0}, c_{0,3},$ and $c_{0,0}$, in terms of $c_{1,3}$, $c_{3,1}$ $c_{2,2}$, $c_{2,3}$, $c_{3,2}, $ and $c_{3,3}$.  (The solutions can be found at the end of the section 4.2.)

Substituting these solutions into our four orthogonality equations gives us four equations and six unknowns.  These, admittedly long, equations are
\begin{eqnarray*}
&& 20c_{3,3}c_{1,3} + 20c_{3,3}c_{3,1} - 4c_{2,2}c_{1,3}-8c_{3,3}c_{2,2}
+ 8c_{3,1}c_{1,3} - 4c_{3,1}c_{2,2} - 8c_{3,3}c_{2,3} - 4c_{3,1}c_{2,3} \\
&& - 4c_{1,3}c_{3,2}
- 8c_{3,3}c_{3,2} 
-4c_{3,1}c_{3,2} - 2c_{1,3} - 2c_{2,2} - 2c_{3,1}
- 2c_{3,3} + 4c_{3,2}^2 + 20c_{3,3}^2 + 8c_{1,3}^2 \\
&& + 4c_{2,2}^2 + 4c_{2,3}^2 + 8c_{3,1}^2 + \frac{5}{2} - 4c_{2,3}c_{1,3}=4
\end{eqnarray*}

\begin{eqnarray*}
&& -c_{3,3}c_{3,1}-c_{3,3}^2+c_{3,3}-c_{3,3}c_{1,3}-c_{3,1}c_{3,2}-2c_{3,3}c_{3,2}
+\frac{1}{2}c_{3,2}-c_{1,3}c_{3,2}+c_{3,2}^2-c_{3,1}c_{2,3} \\
&&-2c_{3,3}c_{2,3}+c_{2,3}^2+
\frac{1}{2}c_{2,3}-c_{2,3}c_{1,3}-c_{3,1}c_{2,2}-2c_{3,3}c_{2,2}+c_{2,2}^2-c_{2,2}c_{1,3}=0 
\end{eqnarray*}

\begin{eqnarray*}
&&-2c_{3,3}c_{1,3}-10c_{3,3}c_{3,1}+2c_{2,2}c_{1,3}+4c_{3,3}c_{2,2}-4c_{3,1}c_{1,3}
+2c_{3,1}c_{2,2}+4c_{3,3}c_{2,3}+2c_{3,1}c_{2,3} \\
&&+2c_{1,3}c_{3,2}+4c_{3,3}c_{3,2}+2c_{3,1}c_{3,2}-c_{1,3}+c_{2,2}+c_{3,1}-c_{3,3}-2c_{3,2}^2-6c_{3,3}^2-2c_{2,2}^2-2c_{2,3}^2-4c_{3,1}^2 \\
&&+2c_{2,3}c_{1,3}+\frac{1}{4}=0 
\end{eqnarray*}
 
\begin{eqnarray*}
&&-10c_{3,3}c_{1,3}-2c_{3,3}c_{3,1}+2c_{2,2}c_{1,3}+4c_{3,3}c_{2,2}-4c_{3,1}c_{1,3}
+2c_{3,1}c_{2,2}+4c_{3,3}c_{2,3}+2c_{3,1}c_{2,3} \\
&&+2c_{1,3}c_{3,2}+4c_{3,3}c_{3,2}+2c_{3,1}c_{3,2}
+c_{1,3}+c_{2,2}-c_{3,1}-c_{3,3}-2c_{3,2}^2-6c_{3,3}^2-4c_{1,3}^2-2c_{2,2}^2-2c_{2,3}^2 \\
&&+2c_{2,3}c_{1,3}+\frac{1}{4}=0 
\end{eqnarray*}

\noindent Notice that all of these equations are quadratic in the six variables, so solutions are basically intersection points of hyper-spheres.

We tried to solve three of the four new equations with \emph{Maple}, but the output contained six possible solution types.  Four of these possible solution types led to refinement coefficients with imaginary components.  The other two gave us real refinement coefficients, but it was clear that we did not discover \emph{all} of the solutions because a solution published in {\bf [H]} was not included.

We then found a better way to solve the fourteen equations.  The method to solve these four equations, giving us the eight types of solutions listed in the next section, is described in detail in the \emph{Maple} file \verb!8solutions.mws!.  Basically, the idea is to eliminate equations by combining them in such a way as to remove quadratic terms.  For instance, subtracting the last two of the four equations above yields this much-simpler equation:

$$-8c_{3,3}c_{1,3}+8c_{3,3}c_{3,1}+2c_{1,3}-2c_{3,1}-4c_{1,3}^2+4c_{3,1}^2 = 0.$$

\noindent  This equation can be solved for $c_{1,3}$, giving us two solutions.  One solution is $c_{1,3}=c_{3,1}$, and we refer to these as type ``A'' solutions.  The other is $c_{1,3} = -2c_{3,3}+\frac{1}{2}-c_{3,1}$, and these are called type ``B'' solutions.  There are four versions of each solution, as described in the next section and proved in the \emph{Maple} file.

\subsection{Method of Creating Masks}

In this section, we describe how to create any possible \emph{mask of refinement coefficients} that satisfy our conditions.  A \emph{mask} is a matrix of values for the 16 refinement coefficients.  To begin, we must recognize that two of the coefficients will become parameters in our final solution, as there are only fourteen equations.  In our solution, the parameters will be $c_{3,2}$ and $c_{3,3}$.

To create a mask, first choose values of $c_{3,2}$ and $c_{3,3}$, along with one of the 8 types of solutions listed below.  Note that there are requirements (described below) in your choice of $c_{3,2}$ and $c_{3,3}$, in order to get real-valued solutions.  Next, compute four of the refinement coefficients in this order:  $c_{2,3}$, $c_{2,2}$, $c_{3,1}$, and then $c_{1,3}$.  Finally, the other refinement coefficients can be computed.

In all of the following solutions, $\mu_1 = \frac{1 + \sqrt{3}}{4}$ and $\mu_2 = \frac{1 - \sqrt{3}}{4}$.  The following discriminants are used:

\begin{eqnarray*}
\Delta_1 &=&  -34+32c_{3,3}\sqrt{3}-20\sqrt{3}+32c_{3,2}\sqrt{3}-48c_{3,3}^2+32c_{3,3}-32c_{3,3}c_{3,2}+48c_{3,2}-48c_{3,2}^2    \\
\Delta_2 &=&  -48c_{3,2}^2-32c_{3,2}\sqrt{3}-32c_{3,3}c_{3,2}+48c_{3,2}-34-32c_{3,3}\sqrt{3}+20\sqrt{3}-48c_{3,3}^2+32c_{3,3}    \\
\Delta_3 &=&  -48c_{3,2}^2-32c_{3,3}c_{3,2}-48c_{3,3}^2+2-16c_{3,3}    \\
\end{eqnarray*}

Here are the eight types of solutions:

Type A1a
\begin{eqnarray*}
c_{2,3} &=& \frac{1}{4}+2\mu_1-\frac{1}{2}(c_{3,2}+c_{3,3})+\frac{1}{8}\sqrt{\Delta_1} \\
c_{2,2} &=& \frac{3}{4}+4\mu_1-c_{3,2}-c_{3,3}-c_{2,3} \\
c_{3,1} &=& \mu_1-c_{3,3}   \\
c_{1,3} &=& c_{3,1}   \\
\end{eqnarray*}

Type A1b
\begin{eqnarray*}
c_{2,3} &=& \frac{1}{4}+2\mu_1-\frac{1}{2}(c_{3,2}+c_{3,3})-\frac{1}{8}\sqrt{\Delta_1} \\
c_{2,2} &=& \frac{3}{4}+4\mu_1-c_{3,2}-c_{3,3}-c_{2,3} \\
c_{3,1} &=& \mu_1-c_{3,3}   \\
c_{1,3} &=& c_{3,1}   \\
\end{eqnarray*}

Requirement on $c_{3,2}$ and $c_{3,3}$ for real solutions of type A1a and A1b:  $\Delta_1 \ge 0$.

Type A2a
\begin{eqnarray*}
c_{2,3} &=& \frac{1}{4}+2\mu_2-\frac{1}{2}(c_{3,2}+c_{3,3})+\frac{1}{8}\sqrt{\Delta_2} \\
c_{2,2} &=& \frac{3}{4}+4\mu_2-c_{3,2}-c_{3,3}-c_{2,3} \\
c_{3,1} &=& \mu_2-c_{3,3}   \\
c_{1,3} &=& c_{3,1}   \\
\end{eqnarray*}

Type A2b
\begin{eqnarray*}
c_{2,3} &=& \frac{1}{4}+2\mu_2-\frac{1}{2}(c_{3,2}+c_{3,3})-\frac{1}{8}\sqrt{\Delta_2} \\
c_{2,2} &=& \frac{3}{4}+4\mu_2-c_{3,2}-c_{3,3}-c_{2,3} \\
c_{3,1} &=& \mu_2-c_{3,3}   \\
c_{1,3} &=& c_{3,1}   \\
\end{eqnarray*}

Requirement on $c_{3,2}$ and $c_{3,3}$ for real solutions of type A2a and A2b:  $\Delta_2 \ge 0$.

Type B1a
\begin{eqnarray*}
c_{2,3} &=& -\frac{1}{2}(c_{3,2}+c_{3,3})+\frac{1}{8}\sqrt{\Delta_3} \\
c_{2,2} &=& \frac{1}{4}-c_{3,2}-c_{3,3}-c_{2,3} \\
c_{3,1} &=& \mu_1-c_{3,3}   \\
c_{1,3} &=& \frac{1}{2}-2c_{3,3}-c_{3,1}   \\
\end{eqnarray*}

Type B1b
\begin{eqnarray*}
c_{2,3} &=& -\frac{1}{2}(c_{3,2}+c_{3,3})-\frac{1}{8}\sqrt{\Delta_3} \\
c_{2,2} &=& \frac{1}{4}-c_{3,2}-c_{3,3}-c_{2,3} \\
c_{3,1} &=& \mu_1-c_{3,3}   \\
c_{1,3} &=& \frac{1}{2}-2c_{3,3}-c_{3,1}   \\
\end{eqnarray*}

Type B1a
\begin{eqnarray*}
c_{2,3} &=& -\frac{1}{2}(c_{3,2}+c_{3,3})+\frac{1}{8}\sqrt{\Delta_3} \\
c_{2,2} &=& \frac{1}{4}-c_{3,2}-c_{3,3}-c_{2,3} \\
c_{3,1} &=& \mu_2-c_{3,3}   \\
c_{1,3} &=& \frac{1}{2}-2c_{3,3}-c_{3,1}   \\
\end{eqnarray*}

Type B1a
\begin{eqnarray*}
c_{2,3} &=& -\frac{1}{2}(c_{3,2}+c_{3,3})-\frac{1}{8}\sqrt{\Delta_3} \\
c_{2,2} &=& \frac{1}{4}-c_{3,2}-c_{3,3}-c_{2,3} \\
c_{3,1} &=& \mu_2-c_{3,3}   \\
c_{1,3} &=& \frac{1}{2}-2c_{3,3}-c_{3,1}   \\
\end{eqnarray*}

Requirement on $c_{3,2}$ and $c_{3,3}$ for real solutions of type B1a, B1b, B2a, and B2b:  $\Delta_3 \ge 0$.  

With these six values, the other ten can then be calculated.  The following equations are the solutions to the 10 linear equations in terms of $c_{1,3}, c_{2,2}, c_{2,3}, c_{3,1}, c_{3,2} \ \textnormal{and} \ c_{3,3}$.

\begin{eqnarray*}
c_{0,0}&=&-c_{3,1}-2c_{3,3}-c_{1,3}+c_{2,2}\\
c_{0,1}&=&-c_{3,1}-2c_{3,3}+\frac{1}{2}+c_{2,3}-c_{1,3}\\
c_{0,2}&=&c_{1,3}+c_{3,3}+\frac{1}{2}-c_{2,2}\\
c_{0,3}&=&c_{3,3}-c_{2,3}+c_{1,3}\\ 
c_{1,0}&=&-c_{3,1}-2c_{3,3}+c_{3,2}-c_{1,3}+\frac{1}{2}\\
c_{1,1}&=&-c_{1,3}-c_{3,1}-c_{3,3}+1\\
c_{1,2}&=&-c_{3,2}+c_{1,3}+c_{3,3}+\frac{1}{2}\\
c_{2,0}&=&c_{3,1}+c_{3,3}-c_{2,2}+\frac{1}{2}\\
c_{2,1}&=&c_{3,1}+c_{3,3}+\frac{1}{2}-c_{2,3}\\
c_{3,0}&=&c_{3,1}+c_{3,3}-c_{3,2}\\
\end{eqnarray*}

An example of a bivariate Daubechies scaling function that can be found in the literature is in {\bf [H]}, a paper published in 1999 by Bin Han of Princeton.  Han's scaling function is of type B2A, using $c_{3,2} = \mu_2$ and $c_{3,3} = 0$.

In the next section, we will prove an important lemma.  Following that is a description of how to generate a picture of a scaling function from a mask.  Then, we will present masks of all eight types, along with the graphs of the related scaling functions.

\section{An Important Lemma}

In this section, we prove another Lemma that will be useful in the rest of this report.

{\bf Lemma 3}: $\sum_i\sum_j \phi(x-i,y-j)=1$.

{\bf Proof of Lemma 3}:  We will first show that
$$\sum_n\sum_m c_{2n,2m}=\sum_n\sum_m c_{2n+1,2m}=\sum_n\sum_m c_{2n+1,2m+1}=\sum_n\sum_m c_{2n,2m+1}=1.$$

Define the following values:
\begin{eqnarray*}
w_1&=&\sum_n\sum_m c_{2n,2m} \\
&=&c_{0,0}+c_{0,2}+c_{2,0}+c_{2,2} \\
w_2&=&\sum_n\sum_m c_{2n+1,2m} \\
&=&c_{1,0}+c_{1,2}+c_{3,0}+c_{3,2} \\
w_3&=&\sum_n\sum_m c_{2n+1,2m+1} \\
&=&c_{1,1}+c_{1,3}+c_{3,1}+c_{3,3} \\
w_4&=&\sum_n\sum_m c_{2n,2m+1} \\
&=&c_{0,1}+c_{0,3}+c_{2,1}+c_{2,3}
\end{eqnarray*}
\noindent The first regularity equation can be written as $w_4-w_1=0$.  Similarly the second and third regularity equations can be written as $w_1-w_4=0$ and $w_2-w_4=0$.  Our averaging equation can be written as $w_1+w_2+w_3+w_4=4$.  We now have a system of four equations with four variables.  The unique solution to this system is $w_1=w_2=w_3=w_4=1$.

Now we will show that $\sum_i\sum_j \phi(x-i,y-j)=1$.  Define
\begin{equation}
f(x,y)=\sum_i\sum_j\phi(x-i,y-j).
\end{equation}
\noindent
Notice that $f(x,y)$ is continuous, as it is defined in terms of the scaling function which is continuous.  Using the dialation equation $f(x,y)$ can also be written as:
\begin{eqnarray*}
f(x,y)&=&\sum_i\sum_j (\sum_m\sum_n c_{m,n}\phi(2(x-i)-m,2(y-j)-n)) \\
&=&\sum_i\sum_j\sum_m\sum_n c_{m,n}\phi(2x-2i-m,2y-2j-n) 
\end{eqnarray*}
Let $a=2i+m$ and $b=2j+n$ so $m=a-2i$ and $n=b-2j$. Thus, 
\begin{eqnarray*}
f(x,y)&=&\sum_i\sum_j\sum_a\sum_b c_{a-2i,b-2j}\phi(2x-a,2y-b) \\
&=&\sum_a\sum_b(\sum_i\sum_j c_{a-2i,b-2j})\phi(2x-a,2y-b)
\end{eqnarray*}
We can see that $\sum_i\sum_j c_{a-2i,b-2j}=1$ regardless of $a$ and $b$ because of our argument above. Thus, 
\begin{eqnarray*}
f(x,y)&=&\sum_a\sum_b\phi(2x-a,2y-b) \\
&=&f(2x,2y)
\end{eqnarray*}

In order to show that $f(x,y)$ is constant, a proof by contradiction will be used.  Assume that $f(0,0)=A$ and for some $p$ and $q$, $f(p,q)=B$ such that $A\not=B$.  Since $f(x,y)=f(2x,2y)$ we can deduce that $f(\frac{p}{2},\frac{q}{2})=B$, $f(\frac{p}{4},\frac{q}{4})=B$, etc.  This pattern continues and since $f(x,y)$ is continuous, it must be that $f(0,0)=B$ which is a contradiction of the assumption that $A\not=B$.  Therefore, $f(x,y)$ must be constant, let's call this constant $P$.

Now we simply need to show that $P$=1.  Since $P$ is a constant, 
\begin{eqnarray*}
P&=&\int_0^1\int_0^1 P dx dy \\
&=&\int_0^1\int_0^1 f(x,y) dx dy 
\end{eqnarray*}
From the definition of $f(x,y)$,
\begin{eqnarray*}
P&=&\int_0^1\int_0^1 \sum_i\sum_j \phi(x-i,y-j) dx dy \\
&=&\sum_i\sum_j \int_0^1\int_0^1 \phi(x-i,y-j) dx dy 
\end{eqnarray*}
The previous step is valid because all but a finite number of terms are zero .  
Using a change of variables,  we now have

\begin{eqnarray*}
P&=&\sum_i\sum_j \int_i^{i+1}\int_j^{j+1} \phi(x,y) dx dy \\
&=&\int_{-\infty}^\infty \int_{-\infty}^\infty \phi(x,y) dx dy\\
&=&1
\end{eqnarray*}
\noindent by the averaging assumption (\ref{E:Averaging2}).  Therefore, we have shown that $\sum_i\sum_j \phi(x-i,y-j)=1$.

Due to the compact support condition for $\phi$, we can actually state that $\phi(2,2)+\phi(1,1)+\phi(2,1)+\phi(1,2)=1$.

\section{Generating $\phi(x,y)$ with the Cascade Algorithm}

\subsection{Overview}

Once a mask is determined by the analysis explained in the previous section, there is a question of how to create the associated scaling function $\phi(x,y)$.  To do so requires a modification of Daubechies' cascade algorithm {\bf [D]}.  For the \emph{bi-cascade algorithm}, we take advantage of the compact support condition, which allows us to state that $\phi(x,y) = 0$ when $x$ and $y$ are integers, except at these four key points:  $(1,1)$, $(1,2)$, $(2,1)$, and $(2,2)$.  We will also make use of the dilation equation (\ref{E:Dilation}).

Basically, there are two parts to the bi-cascade algorithm.  First, determine the values $\phi(1,1)$, $\phi(1,2)$, $\phi(2,1)$ and $\phi(2,2)$ by using a fixed-point method.  Then, determine the values at points $(x,y)$, where $x$ and $y$ are rational numbers whose denominators are powers of 2.  Next, we describe each part in detail, and \emph{Maple} code for the algorithm can be found in the file \verb!bicascade.mws!.

\subsection{Determination of four key function values}

Due to the compact support condition, the dilation equation (\ref{E:Dilation}), evaluated at our four key points, yields the following four equations:
\begin{eqnarray*}
\phi(1,1) &=&  c_{1,1}\phi(1,1)+ c_{0,1}\phi(2,1) + c_{1,0}\phi(1,2) + c_{0,0}\phi(2,2) \\ 
\phi(2,1) &=&  c_{3,1}\phi(1,1)+ c_{2,1}\phi(2,1) + c_{3,0}\phi(1,2) + c_{2,0}\phi(2,2) \\
\phi(1,2) &=&  c_{1,3}\phi(1,1)+ c_{0,3}\phi(2,1) + c_{1,2}\phi(1,2) + c_{0,2}\phi(2,2) \\
\phi(2,2) &=&  c_{3,3}\phi(1,1)+ c_{2,3}\phi(2,1) + c_{3,2}\phi(1,2) + c_{2,2}\phi(2,2) \\
\end{eqnarray*}
\noindent  Another way to look at this is that if $b=[\phi(1,1), \phi(2,1), \phi(1,2), \phi(2,2)]^T$, then we seek solutions of the equation $b=Lb$, where $L$ is the 4-by-4 matrix
$$L = \left (\begin{array}{cccc}
 c_{1,1} & c_{0,1} & c_{1,0} & c_{0,0} \\ c_{3,1} & c_{2,1} & c_{3,0} & c_{2,0} \\ c_{1,3} & c_{0,3} & c_{1,2} & c_{0,2} \\ c_{3,3} & c_{2,3} & c_{3,2} & c_{2,2} \end{array} \right) .$$
\noindent  Notice that any vector $b$ which satisfies the equation is an eigenvector corresponding to eigenvalue 1.

It turns out that one of the eigenvalues of $L$ needs to be 1, while the rest have absolute value less than 1.  Recall that the averaging and regularity equations formed a linear system.  This system determines values for the ten following variables $c_{1,2}, c_{2,1}, c_{2,0}, c_{0,2}, c_{0,1}, c_{1,0}, c_{1,1}, c_{3,0}, c_{0,3},$ and $c_{0,0}$ in terms of $c_{1,3}$, $c_{3,1}$ $c_{2,2}$, $c_{2,3}$, $c_{3,2}, $ and $c_{3,3}$.  Substituting those values into $L$ and computing the eigenvalues with \emph{Maple} yields four eigenvalues:  1, $\frac{1}{2}$ twice, and $-c_{3,2}-c_{2,3}+c_{2,2}+c_{3,3}$.  (See the \emph{Maple} file \verb!eigenL.mws!.)  To show that this last eigenvalue has absolute value less than one, we can use Lagrange multipliers, using the orthogonality equations to get two constraints.  Subject to these constraints, we conclude that
$$|-c_{3,2}-c_{2,3}+c_{2,2}+c_{3,3}| \le \frac{\sqrt{7}}{4} < 1.$$

To determine $b$, we set up an iterative scheme
\begin{equation} \label{E:Iterate}
b_{n+1} = Lb_n.
\end{equation}
\noindent The implication of the analysis of the eigenvalues of $L$ is that this scheme will converege.  Further, in the original cascade algorithm, all function values of $\phi$ on the integers are initially set to zero, except that one of them is set to 1.  This way, the sum of the values of $\phi$ on all of the integers will be 1.  Analogously, for the bi-cascade algorithm, we should initially use $b_0 = [1, 0, 0, 0]^T$, and this is consistent with Lemma 3.

Using \emph{Maple}, the iterative scheme (\ref{E:Iterate}) will converge after about 40 iterations.  We then have the value of $\phi$ at the four key points.

\subsection{Determination of other values of $\phi$}

Once the value of $\phi$ is known at points where $x$ and $y$ are integers, the values where $x$ and $y$ are multiples of $\frac{1}{2}$ can be determined by using the dilation equation (\ref{E:Dilation}).  For example,
$$\phi(\frac{1}{2}, \frac{3}{2}) = \sum_i \sum_j c_{i,j} \phi(1-i, 3-j)$$
\noindent  This process can then be iterated to determine the values of $\phi$ on multiples of $\frac{1}{4}$, $\frac{1}{8}$, etc.  As the iterations continue, the number of points where $\phi$ is determined increases exponentially, and after determining values on multiples of $\frac{1}{64}$, the amount of time needed for the next iterate is not worth the effort.

At the end of the \emph{Maple} file are several methods for creating graphs of the scaling function.  One method is a simple plot of points.  The second method is to define a function which agrees with $\phi$ on all of the determined points.  This function is useful for creating linear combinations of translates of $\phi$, and will be used in a later section.  The third method creates an animated gif that rotates the graph of the scaling functions about the $z$-axis.

\section{Examples of Masks and Associated Scaling Functions}

The following masks were created using the method described in section 4.  Each mask is presented in the following matrix form:

\begin{displaymath}
\left( \begin{array}{cccc}
c_{0,0} & c_{0,1} & c_{0,2} & c_{0,3} \\
c_{1,0} & c_{1,1} & c_{1,2} & c_{1,3} \\
c_{2,0} & c_{2,1} & c_{2,2} & c_{2,3} \\
c_{3,0} & c_{3,1} & c_{3,2} & c_{3,3} \\
\end {array} \right)
\end{displaymath}

In our examples, different parameters were used for the A1 examples than for the other examples.  There is no pair $(c_{3,2}, c_{3,3})$ that will yield real valued masks for all 8 types.  (See Section 4.2.)

A1a, using ($c_{3,2}, c_{3,3}$) =  (1,$\frac{1}{2}$):

{\fontsize{6}{8}

\begin{displaymath} 
\left( \begin{array}{cccc}
-\frac{1}{4} - \frac{1}{8} \sqrt(-46 + 28 \sqrt3) & \frac{1}{8} \sqrt(-46 + 28 \sqrt3) & \frac{1}{2} - \frac{1}{4} \sqrt3 + \frac{1}{8} \sqrt(-46 + 28 \sqrt3)
& \frac{1}{4} - \frac{1}{4} \sqrt3 - \frac{1}{8} \sqrt(-46 + 28 \sqrt3)\\
1 - \frac{1}{2} \sqrt3 & 1 - \frac{1}{2} \sqrt3 & -\frac{1}{4} + \frac{1}{4} \sqrt3 & 
-\frac{1}{4} + \frac{1}{4} \sqrt3(3)\\
\frac{1}{2} - \frac{1}{4} \sqrt3 + \frac{1}{8} \sqrt(-46 + 28 \sqrt3) &  \frac{3}{4} - \frac{1}{4} \sqrt3 - \frac{1}{8} \sqrt(-46 + 28 \sqrt3) & \frac{1}{4} + \frac{1}{2} \sqrt3 - \frac{1}{8} \sqrt(-46 + 28 \sqrt3) & \frac{1}{2} \sqrt3 + \frac{1}{8} \sqrt(-46 + 28 \sqrt3)\\
-\frac{3}{4} + \frac{1}{4} \sqrt3 & -\frac{1}{4} + \frac{1}{4} \sqrt3 & 1 & \frac{1}{2}\\
\end {array} \right)  
\end{displaymath}

}

A1b, using ($c_{3,2}, c_{3,3}$) =  (1,$\frac{1}{2}$):

{\fontsize{6}{8}

\begin{displaymath}
\left( \begin{array}{cccc}
-\frac{1}{4} + \frac{1}{8} \sqrt(-46 + 28 \sqrt3) &  - \frac{1}{8} \sqrt(-46 + 28 \sqrt3)
& \frac{1}{2} - \frac{1}{4} \sqrt3 - \frac{1}{8} \sqrt(-46 + 28 \sqrt3) & \frac{1}{4} - \frac{1}{4} \sqrt3 + \frac{1}{8} \sqrt(-46 + 28 \sqrt3)\\
1 - \frac{1}{2} \sqrt3 & 1 - \frac{1}{2} \sqrt3 & -\frac{1}{4} + \frac{1}{4} \sqrt3 & -\frac{1}{4} + \frac{1}{4} \sqrt3\\
\frac{1}{2} - \frac{1}{4} \sqrt3 - \frac{1}{8} \sqrt(-46 + 28 \sqrt3) & \frac{3}{4} - \frac{1}{4} \sqrt3 + \frac{1}{8} \sqrt(-46 + 28 \sqrt3) & \frac{1}{4} + \frac{1}{2} \sqrt3 + \frac{1}{8} \sqrt(-46 + 28 \sqrt3) & \frac{1}{2} \sqrt3 - \frac{1}{8} \sqrt(-46 + 28 \sqrt3)\\
-\frac{3}{4} + \frac{1}{4} \sqrt3 & -\frac{1}{4} + \frac{1}{4} \sqrt3 & 1 & \frac{1}{2}\\
\end {array} \right)
\end{displaymath}

}

A2a, using ($c_{3,2}, c_{3,3}$) = (0,0):

{\fontsize{6}{8}

\begin{displaymath}
\left( \begin{array}{cccc}
\frac{1}{2} - \frac{1}{8}\sqrt(-34 + 20\sqrt3) & \frac{3}{4} + \frac{1}{8}\sqrt(-34 + 20\sqrt3) &
-\frac{1}{4} + \frac{1}{4}\sqrt3 + \frac{1}{8}\sqrt(-34 + 20\sqrt3) & -\frac{1}{2} + \frac{1}{4}\sqrt3 - \frac{1}{8}\sqrt(-34 + 20\sqrt3) \\
\frac{1}{2}\sqrt3 & \frac{1}{2} + \frac{1}{2}\sqrt3 & \frac{3}{4} - \frac{1}{4}\sqrt3 & \frac{1}{4} - \frac{1}{4}\sqrt3 \\
-\frac{1}{4} + \frac{1}{4}\sqrt3 + \frac{1}{8}\sqrt(-34 + 20\sqrt3) &  \frac{1}{4}\sqrt3 - \frac{1}{8}\sqrt(-34 + 20\sqrt3) & 1 - \frac{1}{2}\sqrt3 - \frac{1}{8}\sqrt(-34 + 20\sqrt3) & 
\frac{3}{4} - \frac{1}{2}\sqrt3 + \frac{1}{8}\sqrt(-34 + 20\sqrt3) \\
\frac{1}{4} - \frac{1}{4}\sqrt3 & \frac{1}{4} - \frac{1}{4}\sqrt3 & 0 & 0 \\
\end {array} \right)
\end{displaymath}

}

A2b, using ($c_{3,2}, c_{3,3}$) = (0,0):

{\fontsize{6}{8}

\begin{displaymath}
\left( \begin{array}{cccc}
\frac{1}{2} + \frac{1}{8}\sqrt(-34 + 20 \sqrt3) & \frac{3}{4} - \frac{1}{8}\sqrt(-34 + 20 \sqrt3) & -\frac{1}{4} + \frac{1}{4}\sqrt3 - \frac{1}{8}\sqrt(-34 + 20 \sqrt3) & -\frac{1}{2} + \frac{1}{4}\sqrt3 + \frac{1}{8}\sqrt(-34 + 20\sqrt3) \\
\frac{1}{2}\sqrt3 &  \frac{1}{2} + \frac{1}{2}\sqrt3 & \frac{3}{4} - \frac{1}{4}\sqrt3 &  \frac{1}{4} - \frac{1}{4} \sqrt3 \\
-\frac{1}{4} + \frac{1}{4} \sqrt3 - \frac{1}{8}\sqrt(-34 + 20 \sqrt3)& \frac{1}{4} \sqrt3 + \frac{1}{8}\sqrt(-34 + 20 \sqrt3) & 1 - \frac{1}{2} \sqrt3 + \frac{1}{8}\sqrt(-34 + 20 \sqrt3)
& \frac{3}{4} - \frac{1}{2} \sqrt3 - \frac{1}{8}\sqrt(-34 + 20 \sqrt3) \\ 
\frac{1}{4} - \frac{1}{4} \sqrt3 & \frac{1}{4} - \frac{1}{4} \sqrt3 & 0 & 0 \\
\end {array} \right)
\end{displaymath}

}

B1a, using ($c_{3,2}, c_{3,3}$) =  (0,0):
\begin{displaymath}
\left( \begin{array}{cccc}
-\frac{1}{4} - \frac{1}{8} \sqrt2 & \frac{1}{8} \sqrt2 & \frac{1}{2} - \frac{1}{4} \sqrt3 + \frac{1}{8} \sqrt2 & \frac{1}{4} - \frac{1}{8} \sqrt2 - \frac{1}{4} \sqrt3 \\
0 & \frac{1}{2} &  \frac{3}{4} - \frac{1}{4} \sqrt3 & \frac{1}{4} - \frac{1}{4} \sqrt3 \\
\frac{1}{2} + \frac{1}{4} \sqrt3 + \frac{1}{8} \sqrt2 &  \frac{3}{4} + \frac{1}{4} \sqrt3 - \frac{1}{8} \sqrt2 & \frac{1}{4} - \frac{1}{8} \sqrt2 & \frac{1}{8} \sqrt2 \\
\frac{1}{4} + \frac{1}{4} \sqrt3 &  \frac{1}{4} + \frac{1}{4} \sqrt3 & 0 & 0 \\
\end {array} \right)
\end{displaymath}

B1b, using ($c_{3,2}, c_{3,3}$) =  (0,0):
\begin{displaymath}
\left( \begin{array}{cccc}
-\frac{1}{4} + \frac{1}{8} \sqrt2 & - \frac{1}{8} \sqrt2 & \frac{1}{2} - \frac{1}{4} \sqrt3 - \frac{1}{8} \sqrt2 & \frac{1}{4} + \frac{1}{8} \sqrt2 - \frac{1}{4} \sqrt3 \\
0 &  \frac{1}{2} &  \frac{3}{4} - \frac{1}{4} \sqrt3 & \frac{1}{4} - \frac{1}{4} \sqrt3\\
\frac{1}{2} + \frac{1}{4} \sqrt3 - \frac{1}{8} \sqrt2 & \frac{3}{4} + \frac{1}{4} \sqrt3 + \frac{1}{8} \sqrt2 & \frac{1}{4} + \frac{1}{8} \sqrt2 & - \frac{1}{8} \sqrt2\\
\frac{1}{4} + \frac{1}{4} \sqrt3 & \frac{1}{4} + \frac{1}{4} \sqrt3 & 0 & 0
\end {array} \right)
\end{displaymath}

B2a, using ($c_{3,2}, c_{3,3}$) =  (0,0):
\begin{displaymath}
\left( \begin{array}{cccc}
-\frac{1}{4} - \frac{1}{8} \sqrt2 & \frac{1}{8} \sqrt2 & \frac{1}{2} + \frac{1}{4} \sqrt3 + \frac{1}{8} \sqrt2 & \frac{1}{4} - \frac{1}{8} \sqrt2 + \frac{1}{4} \sqrt3 \\
0 & \frac{1}{2} & \frac{3}{4} + \frac{1}{4} \sqrt3 & \frac{1}{4} + \frac{1}{4} \sqrt3\\
\frac{1}{2} - \frac{1}{4} \sqrt3 + \frac{1}{8} \sqrt2 & \frac{3}{4} - \frac{1}{4} \sqrt3 - \frac{1}{8} \sqrt2 & \frac{1}{4} - \frac{1}{8} \sqrt2 & \frac{1}{8} \sqrt2\\
\frac{1}{4} - \frac{1}{4} \sqrt3 & \frac{1}{4} - \frac{1}{4} \sqrt3 & 0 & 0\\
\end {array} \right)
\end{displaymath}

B2b, using ($c_{3,2}, c_{3,3}$) =  (0,0):
\begin{displaymath}
\left( \begin{array}{cccc}
-\frac{1}{4} + \frac{1}{8} \sqrt2 & - \frac{1}{8} \sqrt2 & \frac{1}{2} + \frac{1}{4} \sqrt3 - \frac{1}{8} \sqrt2 & \frac{1}{4} + \frac{1}{8} \sqrt2 + \frac{1}{4} \sqrt3 \\
0 & \frac{1}{2} & \frac{3}{4} + \frac{1}{4} \sqrt3 & \frac{1}{4} + \frac{1}{4} \sqrt3\\
\frac{1}{2} - \frac{1}{4} \sqrt3 - \frac{1}{8} \sqrt2 & \frac{3}{4} - \frac{1}{4} \sqrt3 + \frac{1}{8} \sqrt2 & \frac{1}{4} + \frac{1}{8} \sqrt2 & - \frac{1}{8} \sqrt2 \\
\frac{1}{4} - \frac{1}{4} \sqrt3 & \frac{1}{4} - \frac{1}{4} \sqrt3 & 0 & 0\\
\end {array} \right)
\end{displaymath}

\section{Reproducing Linear Functions}

\subsection{A General Formula to Reproduce Linear Functions}

Just like $D_4$, we can reproduce linear functions by using linear combinations of $\phi(x,y)$ and its translates.  To begin, suppose $g$ is a function that we wish to reproduce, that is, 
$$g(x,y)=\sum_u\sum_va_{u,v}\phi(x-u,y-v)$$
\noindent With a change of variables, we can rewrite this equation as
\begin{equation} \label{E:LinComb}
g(x,y)=\sum_i\sum_ja_{x-i,y-j}\phi(i,j)      
\end{equation}
So, using integer values for $x$ and $y$ and an arbitrary bivariate $\phi$ mask,
\begin{equation} \label{E:Sumof4}
g(x,y)=a_{x-2,y-2}\phi(2,2)+a_{x-1,y-1}\phi(1,1) + a_{x-2,y-1}\phi(2,1)+a_{x-1,y-2}\phi(1,2)
\end{equation}
Notice that we only have $\phi(2,2), \phi(2,1), \phi(1,2), $and $\phi(1,1)$ left because $\phi(i,j)=0$ at all other integer points.

We are now ready to figure out $a_{x-i,y-j}$ in terms of $\phi(i,j)$.  This will enable us to create linear functions out of a bivariate scaling function $\phi$ and its translates.  We assume that $g(x,y)$ is a polynomial, so $a_{u,v}$ as a function of $u$ and $v$ will have the same degree as $g(x,y)$.  

If $g(x,y)=x$, then $a_{i,j}=\alpha_xi+\beta_xj+\gamma_x$ for some constants $\alpha_x, \beta_x,$ and $\gamma_x$. We will determine these constants such that we can write $x$ in terms of $\phi(2,2),\phi(1,1),\phi(2,1)$ and $\phi(1,2)$. Substituting into (\ref{E:Sumof4}), and using Lemma 3:
\begin{eqnarray*}
x &=&a_{x-2,y-2}\phi(2,2)+a_{x-1,y-1}\phi(1,1)+a_{x-2,y-1}\phi(2,1)+a_{x-1,y-2}\phi(1,2)\\
&=&(\alpha_x(x-2)+\beta_x(y-2)+\gamma_x)\phi(2,2)+(\alpha_x(x-1)+\beta_x(y-1)+\gamma_x)\phi(1,1)\\
&&+(\alpha_x(x-2)+\beta_x(y-1)+\gamma_x)\phi(2,1)+(\alpha_x(x-1)+\beta_x(y-2)+\gamma_x)\phi(1,2) \\
&=&x\alpha_x(\phi(2,2)+\phi(1,1)+\phi(2,1)+\phi(1,2))+y\beta_x(\phi(2,2)+\phi(1,1)+\phi(2,1)+\phi(1,2))\\
&&+\alpha_x(-2\phi(2,2)-\phi(1,1)-2\phi(2,1)-\phi(1,2))+\beta_x(-2\phi(2,2)-\phi(1,1)-\phi(2,1)-2\phi(1,2))+\\
&&\gamma_x(\phi(2,2)+\phi(1,1)+\phi(2,1)+\phi(1,2)) \\
&=&x\alpha_x+y\beta_x+\alpha_x(-2\phi(2,2)-\phi(1,1)-2\phi(2,1)-\phi(1,2)) \\
&&+\beta_x(-2\phi(2,2)-\phi(1,1)-\phi(2,1)-2\phi(1,2))+\gamma_x \\
\end{eqnarray*}
\noindent Therefore $\alpha_x=1$ and $\beta_x=0$, and $\gamma_x=2\phi(2,2)+\phi(1,1)+2\phi(2,1)+\phi(1,2)$. Using these constants we get

\begin{equation} \label{E:FormForX}
x=\sum_i\sum_j(i+2\phi(2,2)+\phi(1,1)+2\phi(2,1)+\phi(1,2))\phi(x-i,y-j)
\end{equation}

Similarly, we will write $y$ in terms of $\phi(2,2),\phi(1,1),\phi(2,1)$ and $\phi(1,2)$.  If $g(x,y)=y$, then $a_{i,j}=\alpha_yi+\beta_yj+\gamma_y$  for some constants $\alpha_y, \beta_y,$ and $\gamma_y$.  Again using (\ref{E:Sumof4}) and Lemma 3 gives us:
\begin{eqnarray*}
y &=&a_{x-2,y-2}\phi(2,2)+a_{x-1,y-1}\phi(1,1)+a_{x-2,y-1}\phi(2,1)+a_{x-1,y-2}\phi(1,2)\\
&=&(\alpha_y(x-2)+\beta_y(y-2)+\gamma_y)\phi(2,2)+(\alpha_y(x-1)+\beta_y(y-1)+\gamma_y)\phi(1,1)\\
&&+(\alpha_y(x-2)+\beta_y(y-1)+\gamma_y)\phi(2,1)+(\alpha_y(x-1)+\beta_y(y-2)+\gamma_y)\phi(1,2) \\
&=&x\alpha_y(\phi(2,2)+\phi(1,1)+\phi(2,1)+\phi(1,2))+y\beta_y(\phi(2,2)+\phi(1,1)+\phi(2,1)+\phi(1,2))\\
&&+\alpha_y(-2\phi(2,2)-\phi(1,1)-2\phi(2,1)-\phi(1,2))+\beta_y(-2\phi(2,2)-\phi(1,1)-\phi(2,1)-2\phi(1,2))+\\
&&\gamma_y(\phi(2,2)+\phi(1,1)+\phi(2,1)+\phi(1,2)) \\
&=&x\alpha_y+y\beta_y+\alpha_y(-2\phi-\phi(1,1)-2\phi(2,1)-\phi(1,2))+ \\
&&\beta_y(-2\phi-\phi(1,1)-2\phi(2,1)-\phi(1,2))+\gamma_y 
\end{eqnarray*}
\noindent Therefore $\alpha_y=0$, $\beta_y=1$, and $\gamma_y = 2\phi(2,2)+\phi(1,1)+\phi(2,1)+2\phi(1,2)$.  Thus, we have found constants $\alpha_y,\beta_y,$ and $\gamma_y$ such that $g(x,y)=y$, and
\begin{equation} \label{E:FormForY}
y=\sum_i\sum_j(j+2\phi(2,2)+\phi(1,1)+\phi(2,1)+2\phi(1,2))\phi(x-i,y-j)
\end{equation}

In equations (\ref{E:FormForX}) and (\ref{E:FormForY}), along with Lemma 3, we now have $x$, $y$ and $1$ written in terms of $\phi(2,2),\phi(1,1),\phi(2,1),$ and $\phi(1,2)$ .  This enables us to write an approximation for any type of plane given a particular mask.

Let us look at the general case of $f(x,y)=kx+ly+m$ where $k,l,$ and $m$ are constants.  Thus,
\begin{eqnarray*}
f(x,y)&=&kx+ly+m\\
&=&k(\sum_i\sum_j(i+2\phi(2,2)+\phi(1,1)+2\phi(2,1)+\phi(1,2))\phi(x-i,y-j)) \\
&&+l(\sum_i\sum_j(j+2\phi(2,2)+\phi(1,1)+\phi(2,1)+2\phi(1,2))\phi(x-i,y-j)) \\
&&+m(\sum_i\sum_j\phi(x-i,y-j)) \\
&=&\sum_i\sum_j (k(i+2\phi(2,2)+\phi(1,1)+2\phi(2,1)+\phi(1,2))+l(j+2\phi(2,2)+\phi(1,1)+ \\
&&\phi(2,1)+2\phi(1,2))+m)\phi(x-i,y-j) \\
\end{eqnarray*}
\noindent  This simplifies to
\begin{equation} \label{E:ReproduceG}
kx+ly+m = \sum_i\sum_j (2\phi(2,2)(k+l)+\phi(1,1)(k+l)+\phi(2,1)(2k+l)+\phi(1,2)(k+2l) \\
+ki+lj+m)\phi(x-i,y-j)
\end{equation}
\noindent Using the above equation, we can generate $f(x,y)$ using a bivariate scaling function associated with a given mask.

\subsection{An Example of Reproducing A Linear Function}

Consider the function $g(x,y)=x+6y-10$, a linear function that we wish to reproduce using one of our scaling functions and its translates.  In this case, $k=1, l=6,$ and $m=-10$. So,
\begin{eqnarray*}
g(x,y)&=&\sum_i\sum_j (2(1)+2(6))\phi(2,2)+((1)+(6))\phi(1,1)+(2(1)+(6))\phi(2,1) \\ &&+((1)+2(6))\phi(1,2)+(1)i+(6)j+(-10))\phi(x-i,y-j) \\
&=&\sum_i\sum_j (14\phi(2,2)+7\phi(1,1)+8\phi(2,1)+13\phi(1,2)+i+6j-10)\phi(x-i,y-j)
\end{eqnarray*}

We will reproduce $g$ using the scaling function in the figure below.  This is a B2b1 type of scaling function, when $c_{3,2}=c_{3,3}=0$.  Note that this scaling function is also plotted in section 4.2.

$$
\left( \begin{array}{cccc}

-\frac{1}{4}-\frac{1}{8}\sqrt{2} & \frac{1}{8}\sqrt{2} & \frac{1}{2}+\frac{1}{4}\sqrt{3}+\frac{1}{8}\sqrt{2}& 
\frac{1}{4}-\frac{1}{8}\sqrt{2}+\frac{1}{4}\sqrt{3} \\

0 & \frac{1}{2} & \frac{3}{4}+\frac{1}{4}\sqrt{3} & \frac{1}{4}+\frac{1}{4}\sqrt{3} \\

\frac{1}{2}-\frac{1}{4}\sqrt{3}+\frac{1}{8}\sqrt{2} & 
\frac{3}{4}-\frac{1}{4}\sqrt{3}-\frac{1}{8}\sqrt{2} & 
\frac{1}{4}-\frac{1}{8}\sqrt{2} & \frac{1}{8}\sqrt{2} \\

\frac{1}{4}-\frac{1}{4}\sqrt{3} & \frac{1}{4}-\frac{1}{4}\sqrt{3} & 0 & 0 \\

\end {array} \right)
$$

\section{Areas for Further Study}

One area for further study is exploring the mother wavelets.  (Once a scaling function is defined, a mother wavelet can also be defined.  See {\bf [AS]} and {\bf [D]}.)  According to Daubechies, there should be at least three for each of our bivariate scaling function $\phi$.  Another idea is to develop a further classification scheme for these scaling functions.  A bivariate version of $D_6$ would also be interesting, as they could be used to exactly reproduce bowl and saddle surfaces.  Also, our undocumented investigations into bivariate Haar wavelets can be tied into the work in this report.

\section{References}

\noindent {\bf [AS]}  Aboufadel, E. and Schlicker, S., \emph{Discovering Wavelets}, 1999, Wiley.

\noindent {\bf [D]}  Daubechies, I., \emph{Ten Lectures on Wavelets}, 1992, SIAM.

\noindent {\bf [H]}  Han, B., ``Analysis and Construction of Optimal Multivariate Biorthogonal Wavelets with Compact Support'', \emph{SIAM J. Math. Anal.}, 1999, Vol. 31, p.274--304.

\noindent {\bf [J]}  Jia, R.Q., ``Approximation Properties of Multivariate Wavelets'', \emph{Math. Comp.}, 1998, Vol. 67, p. 647--665.

\section{Note about Maple files and Image Files}

A longer version of this paper is available which includes image files of surface graphs of bivariate scaling functions.  It can be downloaded at:

\verb!www.gvsu.edu/mathstat/wavelets/undergrad.htm!

Maple files may be requested from Edward Aboufadel.

\end{document}